\documentclass[a4paper,10pt]{article}

\usepackage{amssymb}
\usepackage{amsmath}
\usepackage{graphicx}
\usepackage{xcolor}
\usepackage{natbib}
\usepackage{algorithm}
\usepackage{algpseudocode}
\usepackage{algorithmicx}
\usepackage[left=2cm,top=2.5cm,right=2cm,bottom=2.5cm]{geometry}

\graphicspath{{Figures/}}

\title{Algorithmic analysis towards time-domain extended source waveform inversion\thanks{Published by \emph{Pure and Applied Geophysics, doi:10.1007/s00024-024-03556-3}}}
\author{Pengliang Yang$^1$ and Wei Zhou$^2$\\
  $^1$ Harbin Institute of Technology, China. E-mail: ypl.2100@gmail.com\\
$^2$ King Fahd University of Petroleum and Minerals, Saudi Arabia. E-mail: wei.zhou@kfupm.edu.sa}

\begin{document}

\maketitle

\begin{abstract}Full waveform inversion (FWI) updates the subsurface model from an initial model by comparing observed and synthetic seismograms. Due to high nonlinearity, FWI is easy to be trapped into local minima. Extended domain FWI, including wavefield reconstruction inversion (WRI) and extended source waveform inversion (ESI) are attractive options to mitigate this issue. 
  This paper makes an in-depth analysis for FWI in the extended domain, identifying key challenges and searching for potential remedies towards practical applications. WRI and ESI are formulated within the same mathematical framework using Lagrangian-based adjoint-state method with a special focus on time-domain formulation using extended sources, while putting connections between classical FWI, WRI and ESI: both WRI and ESI can be viewed as weighted versions of classic FWI. Due to symmetric positive definite Hessian, the conjugate gradient is explored to efficiently solve the normal equation in a matrix free manner, while both time and frequency domain wave equation solvers are feasible. This study finds that the most significant challenge comes from the huge storage demand to store time-domain wavefields through iterations. To resolve this challenge, two possible workaround strategies can be considered, i.e., by extracting sparse frequencial wavefields or by considering time-domain data instead of wavefields for reducing such challenge. We suggest that these options should be explored more intensively for tractable workflows.
\end{abstract}

\section{Introduction}

Full waveform inversion (FWI), proposed by \citet{Lailly_1983_SIP} and \citet{Tarantola_1984_ISR,Tarantola_1984_LIS}, is a high-resolution seismic imaging technique that iteratively minimizes the misfit between the observed seismic data and the synthetic data modelled from the wave equation. Due to its high non-linearity, successful FWI often requires a good starting model with accurate background information, high-quality data with sufficient low-frequency content \citep{Bunks_1995_MSW}. This is partially because the searching space of this inverse problem is reduced to the model space only, while the wavefield is the numerically exact solution to the wave equation. A consequence of this reduced formulation of inversion is the local minima issue that presents FWI from converging to the global minimum. One typical phenomenon coming from the local minima issue is cycle skipping, that the phases between observed and synthetic seismograms are larger than half of the period \citep{Virieux_2009_OFW}.

In order to overcome these issues, significant effort has been performed with the aim to enlarge the basin of the attraction of the misfit functional, while mitigating the strong dependence of the inversion to an accurate initial model. In order to obtain a reliable initial macro-model, \citet{Luo_1991_WETc} developed the wave-equation travel-time inversion using the cross-correlation between the observed and synthetic data to measure and decrease the time shifts of the events. \citet{Bunks_1995_MSW} highlighted the importance of low-frequency content in driving FWI out of local minima, therefore suggesting a multiscale inversion strategy which has become a common practice of the workflow. \citet{Shin_2008_WIL} developed the Laplace domain FWI which highlights the early arrivals and damps the late reflections in retrieving a reliable macro model \citep{Shin_2008_CBB}.

Intensive work by different research groups has been done by modifying the misfit functional to enlarge the width of the convexity valley. \citet{Bozdag_2011_MFF} have studied FWI using instantaneous phase and envelope measurements. \citet{Hale_2013_DWS,Ma_2013_WRT} proposed dynamic warping to capture the same events in the observed and synthetic data and measure their traveltime differences. A deconvolution-based misfit functional for FWI, coined as adaptive waveform inversion, has been proposed by \cite{Warner_2014_AWI} to mitigate the cycle-skipping issue, which has been shown to be equivalent to the introduction of a matching filter \citep{Zhu_2016_BSM}. \citet{Wu_2014_SEI} showed that the envelope-based misfit functional is also a potential candidate for retrieving reliable velocity models in tomographic mode, which could be further improved by cascading a classic FWI.
Optimal transport distance \citep{Engquist_2011_SPH,Metivier_2016_TOF} measures the data difference in a more global manner, allowing translation along temporal and spatial axes of the data to be compared. This distance, sometimes called Earth Mover's Distance (EMD), makes the FWI workflow less sensitive to the initial model design.

Another popular approach for overcoming local minima in linearized optimization is an arbitrary extension of the searching space, from the model space only to a model-wavefield joint space. \citet{VanLeeuwen_2013_MLM} introduced a wavefield reconstruction inversion (WRI) method, and showed that it is more immune to the cycle skipping issue. However, \citet{aghamiry2021efficient} pointed that WRI is mathematically equivalent to the contrast source inversion \citep{abubakar2009application}, for which the local minima issue and cycle skipping persist \citep{van2001contrast}. The misfit functional of WRI includes the wave equation as a penalty. This relaxes the requirement that the wave equation has to be exactly satisfied in each model update iteration. Unlike classical FWI, the objective of WRI is a two-variable optimization including the physical parameters and the seismic wavefield. The penalty formulation leads to a normal equation that involves a Hessian matrix regularized by the penalty parameter. In 2D geometry, such an equation could be solved by direct or iterative frequency-domain solvers, thanks to the light memory requirement for storing the wavefield of few frequency components. However, when a large number of frequencies are required, or when 3D geometry should be considered, frequency-domain solvers become less sufficient than time stepping solvers. Therefore, time-domain implementations of WRI or any extension methods are still interesting for real applications.
Unfortunately, solving the normal equation in the time domain requires saving all snapshots of the wavefield, leading to a huge memory demand \citep{aghamiry2020accurate}. This has become a bottleneck in applying these methods to real problems.

Following the space extension strategy and the differential semblance optimization (DSO) to achieve global convergence \citep{Symes_1998_ASP}, \citet{Symes_2008_MVA,symes2015algorithmic} developed the so-called extended domain FWI (ESI). The key idea is to add extra degrees of freedom to the searching space, which may be unphysical, but effectively fit the observed data through modelling. Then, the optimization aims to decrease these added freedom through iterations without compromising the data fit. Various extension strategies, including source extension \citep{Huang_2016_MSW}, offset extension \citep{fu2017discrepancy}, source and receiver extension \citep{Huang_2017_FWI} and many others \citep{Biondi_2013_TFI,barnier2020full}, have been shown under different settings that the local minima issue of FWI can be mitigated.

Attempts have been performed to connect classical FWI, WRI and ESI \citep{van2019note,symes2020wavefield}. In this paper, we will formulate WRI and ESI in a unified framework, using the Lagrangian-based adjoint-state method. This framework sheds a new light and corrects some misunderstanding among WRI, ESI, and FWI approaches: the misfit functional of both WRI and ESI can be expressed as the weighted  misfit functional of FWI where the weighting matrix is parameter dependent.
We will use a conjugate gradient method to implement time-domain WRI and ESI. We will aso discuss on computationally efficient and memory-frugal methods to fulfills the low memory promises. This study aims to steer the focus of the community from algorithmic issue to the development of efficient low-memory extended domain methods.

\section{FWI, WRI, and ESI formulations revisited}

In this section, the three methods will be reviewed using same notations. Their similarity and differences will be highlighted.

\subsection{Least-squares Full Waveform Inversion}

Classical FWI in a least-squares sense can be expressed by minimizing the following misfit functional with the wave-equation constraint:
\begin{equation}\label{eq:fwi}
  J_\text{FWI}(m) = \frac{1}{2} \|Ru -d\|^2, \qquad \text{s.t.} \quad A(m)u = f
\end{equation}
where the quantity $d:=d(x_r,t;x_s)$ is the observed data recorded at the receiver positions $x_r$, due to a source at the position $x_s$. $u:=u(x,t;x_s)$ is the modeled wavefield. $x\in X$, $t\in[0,T]$. The projection operator $R$ samples $u$ at the receiver location, i.e, $Ru:=\int_X u(x,t;x_s)\delta(x-x_r)dx=u(x_r,t;x_s)$. $A(m)$ is the wave-equation operator with $m$ the target model parameter to be inverted, and $f(t)$ is the source time function. Considering the 2nd-order acoustic equation with a constant density
\begin{equation} \label{eq:wave}
  A(m)u = \left(\frac{1}{v^2} \partial_t^2  -\Delta\right) u = f,
\end{equation}
the model parameter $m$ can be the P-wave velocity $v$.
The minimization of $J_\text{FWI}(m)$ in \eqref{eq:fwi} can be translated into an unconstrained optimization through the adjoint state method based on the Lagrangian formulation \citep{Plessix_2006_RAS}. The Lagrangian function for $J_\text{FWI}(m)$ is given by
\begin{equation}
  \mathcal{L}(u,\lambda,m) = \frac{1}{2}\|Ru -d\|^2 +\langle \lambda,A(m)u- f\rangle_{X\times[0,T]},
\end{equation}
where the Lagrangian multiplier $\lambda=\lambda(x,t)$ is also called the adjoint wavefield. The inner product $\langle \lambda, u\rangle_\Omega$, for any set $\Omega$, is defined as
\begin{equation}
\langle \lambda, u\rangle_\Omega = \int_\Omega \lambda\, u\, \mathrm{d}\Omega.
\end{equation}
The gradient of the Lagrangian is
\begin{equation}
\frac{\mathrm{d} \mathcal{L}}{\mathrm{d} m} = \frac{\partial\mathcal{L}}{\partial u}\frac{\partial u}{\partial m} + \frac{\partial\mathcal{L}}{\partial \lambda}\frac{\partial\lambda}{\partial m} + \langle \lambda, \frac{\partial A}{\partial m}  u\rangle_{[0,T]}.
\end{equation}
Zeroing $\frac{\partial \mathcal{L}}{\partial \lambda}$ recovers the wave equation \eqref{eq:wave}.
Zeroing $\frac{\partial \mathcal{L}}{\partial u}$ leads to the adjoint equation:
\begin{equation} \label{eq:FWI_lambda}
  A^H \lambda = -R^H(Ru-d),
\end{equation}
where $^H$ denotes the adjoint (complex conjugate) operator.
At the saddle point, the gradient of $J_\text{FWI}$ with respect to $m$ is given by:
\begin{equation}\label{eq:gradient}
  \nabla_m J_\text{FWI} = \frac{\partial \mathcal{L}}{\partial m} = \langle \lambda, \frac{\partial A}{\partial m}  u\rangle_{[0,T]} = \int_0 ^T \lambda(x,t) \partial_t^2 u(x,t) \mathrm{d} t.
\end{equation}

\subsection{Wavefield Reconstruction Inversion}
\cite{VanLeeuwen_2013_MLM} takes the wave equation \eqref{eq:wave} as a penalty to regularize the least-squares data error, which leads to the following WRI misfit functional
\begin{equation}
  J_\text{WRI}(m) =\frac{1}{2}\|R u - d\|^2 + \frac{\beta}{2} \|A(m)u -f \|^2 
                  =\frac{1}{2}\left\|\begin{pmatrix}
      R \\
      \sqrt{\beta} A(m)
      \end{pmatrix}u -\begin{pmatrix}
      d\\
      \sqrt{\beta} f
      \end{pmatrix}
                  \right\|^2,
\end{equation}
where both model parameter $m$ and the wavefield $u$ are considered as the parameter of the optimization, while $\beta$ controls the penalty from the wave equation. During the iterations, $\beta$ should be increased so that the wavefield $u$ converges to the solution of the wave equation, see more details in section \ref{sec:beta}. A variable projection by zeroing $\partial J_\text{WRI}/\partial u=0$ gives the normal equation \citep{VanLeeuwen_2013_MLM}:
\begin{equation}\label{eq:normal}
  ( \beta A^H A + R^H R)u = R^H d +  \beta A^H f,
\end{equation}
which permits easy solutions of $u$ in the frequency domain using direct or iterative solvers. However, it is not straightforward to express the operator $(\beta A^H A + R^H R)$ in the time domain, thus preventing the use of time-stepping methods to solve for $u$.

To find a time-domain expression using the adjoint state method, let us introduce a new quantity $q(x,t;x_s)$ to represent the mismatch of the wave equation, that is,
\begin{equation}\label{eq:newq}
 q := A(m)u -f.
\end{equation}
Then, $J_\text{WRI}$ can be rewritten as
\begin{equation}\label{eq:wrimisfit}
J_\text{WRI}(m)=\frac{1}{2}\|R u - d\|^2 + \frac{\beta}{2} \|q\|^2.
\end{equation}
The associated Lagrangian is
\begin{equation}
  \mathcal{L}(m,u,q,\lambda) = \frac{1}{2}\|Ru -d\|^2 +\frac{\beta}{2} \|q\|^2 +\langle \lambda,A(m)u- f-q\rangle_{X\times[0,T]},
\end{equation}
and
\begin{equation}
\frac{\partial \mathcal{L}}{\partial m} =\frac{\partial\mathcal{L}}{\partial u}\frac{\partial u}{\partial m}
  + \frac{\partial\mathcal{L}}{\partial q}\frac{\partial q}{\partial m}
  + \frac{\partial\mathcal{L}}{\partial \lambda}\frac{\partial \lambda}{\partial m}
  + \langle \lambda, \frac{\partial A}{\partial m}  u\rangle_{[0,T]}.
\end{equation}
Zeroing the partial derivatives, we have
\begin{eqnarray}
  \frac{\partial\mathcal{L}}{\partial u}      =0 &\Rightarrow & R^H(Ru-d) + A^H \lambda = 0; \label{eq:newadj1}\\
  \frac{\partial\mathcal{L}}{\partial q}      =0 &\Rightarrow & \beta q -\lambda =0; \label{eq:newadj2} \\
  \frac{\partial\mathcal{L}}{\partial \lambda}=0 &\Rightarrow & A(m)u-f-q=0 \quad \text{i.e. Eq \eqref{eq:newq}.} \label{eq:newadj3}
\end{eqnarray}
At the saddle point, the gradient of $J_\text{WRI}$ with respect to model parameters $m$ is given by
\begin{equation}\label{eq:gradient}
  \nabla_m J_\text{WRI} = \frac{\partial \mathcal{L}}{\partial m} = \langle \lambda, \frac{\partial A}{\partial m}  u\rangle_{[0,T]}.
\end{equation}
Note that, because the wave equation is not strictly satisfied, $u$ and $\lambda$ as given here are different from those quantities in classical FWI (\eqref{eq:wave} and \eqref{eq:FWI_lambda}).

Our formulation is consistent with \cite{VanLeeuwen_2013_MLM}. For example, by \eqref{eq:newadj2} and \eqref{eq:newq}, we have
\begin{equation}\label{eq:ulambda}
A u = f + \frac{1}{\beta}\lambda.
\end{equation}
Applying $A^H$ on both sides and substituting for $\lambda$ in \eqref{eq:newadj1} recovers the normal equation \eqref{eq:normal} of \cite{VanLeeuwen_2013_MLM}. Using \eqref{eq:ulambda}, we can rewrite our gradient expression \eqref{eq:gradient} as
\begin{equation}\label{eq:gradient2}
\nabla_m J_\text{WRI} = \beta\langle Au -f, \frac{\partial A}{\partial m}  u\rangle_{[0,T]},
\end{equation}
which is equivalent to $\beta\sum_\omega (Au -f)^H \frac{\partial A}{\partial m} u$ in the frequency domain. Particularly, when we specify $\frac{\partial A}{\partial m}=-\omega^2$, the expression is same as the one given by \citet[their equation 7]{VanLeeuwen_2013_MLM}, except the scaling constant $\beta$.

A link between $J_\text{WRI}$ and $J_\text{FWI}$ of equation \eqref{eq:fwi} can be shown by defining the Green's function operator $S:=RA^{-1}$.
From equations \eqref{eq:newadj1} and \eqref{eq:newadj2}, we obtain
\begin{equation}
  q = \frac{1}{\beta}\lambda = \frac{1}{\beta}A^{-H} R^H(d-Ru) =\frac{1}{\beta} S^H \Delta d,
\end{equation}
where $\Delta d:=d-Ru$ denotes the data residual. Using this notation, $J_\text{WRI}$ becomes
\begin{equation}\label{eq:wrieqfwi}
    J_\text{WRI}=\frac{1}{2}\|\Delta d\|^2 + \frac{1}{2\beta} \|S^H \Delta d\|^2
    =\frac{1}{2} \Delta d^H  \underbrace{(I + \frac{1}{\beta} SS^H)}_{W_1}\Delta d
    =\frac{1}{2} \|\Delta d\|_{W_1}^2.
\end{equation}
Comparing with $J_\text{FWI}$, it shows that the WRI approach is a least-squares FWI based on weighted data residuals (\citet{van2019note} and \citet[equation 22]{symes2020wavefield}).

Another link between WRI and FWI can be shown by rewriting $J_\text{WRI}$ using the $S$ operator:
\begin{equation}\label{eq:wridSf}
    J_\text{WRI} = \frac{1}{2}\| Sq-(d-Sf) \|^2 + \frac{\beta}{2} \|q\|^2,\\
\end{equation}
which leads to the normal equation for $q$ with an extended source on the right-hand side:
\begin{equation}\label{eq:normalwri}
(S^H S + \beta I) q = S^H (d-S f).
\end{equation}
Let us assume that the initial background model allows to reproduce the first arrivals (e.g. diving waves) with the physical source $f$. Then, the remaining data in $d-Sf=d-RA^{-1}f$ are mainly reflections. In such a scenario, minimizing $J_\text{WRI}$ as expressed in \eqref{eq:wridSf} is to fit the reflection data with extended searching space (presence of $q$), which can be interpreted as reflection waveform inversion based on the penalty formulation.

The above interpretations show that WRI is a specific form of FWI under the Lagrangian formulations. The WRI gradient can still be constructed by the zero-lag crosscorrelation between a forward field and an adjoint field. Since WRI is equivalent to classic FWI by weighting the data residual using $W_1$, it naturally cycle-skips, sufferring the local minima issue (see the 1D example as given by \citet{symes2020wavefield}), despite positive stories \citep{VanLeeuwen_2013_MLM,peters2014wave,aghamiry2020accurate}.

\subsection{Extended Source Inversion}

The key idea of extended source waveform inversion (ESI) proposed by \citet{symes2015algorithmic} is slightly different from the one embedded into the WRI. An extended source $q(x,t;x_s)$ is introduced such that
\begin{equation}\label{eq:newq2}
q:= Au,
\end{equation}
which is almost identical to equation \eqref{eq:newq} in WRI except that the physical source $f$ is absent. This leads to a source-independent inversion scheme.  The corresponding misfit functional is defined as \citep{huang2018source} 
\begin{equation}\label{eq:esi}
J_\text{ESI}(m) = \frac{1}{2} \|Ru-d\|^2 + \frac{\beta}{2}\|Bq\|^2,
\end{equation}
where the penalty term is different from the one in the WRI misfit functional \eqref{eq:wrimisfit} as it involves a weighting operator $B$, also known as the annihilator \citet{Symes_2008_MVA}. Possible options for $B$, e.g. $B=\|x-x_s\|$ and $B=t$ in the form of pseudo-differential operator, can be considered to achieve good convergence \citep{symes2020full}. The associated Lagrangian is
\begin{equation}
\mathcal{L}(u,q,\lambda,m) = \frac{1}{2} \|Ru-d\|^2 + \frac{\beta}{2}\|Bq\|^2 + \langle\lambda, Au -q \rangle_{X\times[0,T]}.
\end{equation}
Zeroing the partial derivatives, we have
\begin{eqnarray}
  \frac{\partial\mathcal{L}}{\partial u}      =0 &\Rightarrow & R^H(Ru-d) + A^H \lambda = 0; \label{eq:ewiadj}\\
  \frac{\partial\mathcal{L}}{\partial q}      =0 &\Rightarrow & \beta B^HB q -\lambda =0; \label{eq:ewilq} \\
  \frac{\partial\mathcal{L}}{\partial\lambda}=0 &\Rightarrow & A(m)u-q=0 \quad \text{i.e. Eq \eqref{eq:newq2}.} \label{eq:ewifwd}
\end{eqnarray}
At the saddle point, the gradient of $J_\text{ESI}$ with respect to model parameters $m$ is given by
\begin{equation}\label{eq:gradewi}
  \nabla_m J_\text{ESI} = \frac{\partial \mathcal{L}}{\partial m} = \langle \lambda, \frac{\partial A}{\partial m}  u\rangle_{[0,T]}.
\end{equation}
Note that, $u$ and $\lambda$ as given here are different from those quantities in classical FWI (\eqref{eq:wave} and \eqref{eq:FWI_lambda}) and WRI (\eqref{eq:newadj1}--\eqref{eq:newadj3}).

By \eqref{eq:ewiadj} and \eqref{eq:ewilq}, we can express the extended source as
\begin{equation}
q=(\beta B^H B)^{-1}\lambda = (\beta B^H B)^{-1} A^{-H} R^H(d-Ru)=(\beta B^H B)^{-1} S^H \Delta d,
\end{equation}
where $\beta B^H B$ is assumed to be invertible. It implies that the definition \eqref{eq:esi} is equivalent to
\begin{equation}
  J_\text{ESI} = \frac{1}{2}\|\Delta d\|^2 + \frac{\beta}{2}\|B(\beta B^H B)^{-1} S^H \Delta d\|^2
  =\frac{1}{2} \Delta d^H \underbrace{(I + \frac{1}{\beta} S(B^H B)^{-1}S^H)}_{W_2} \Delta d
  =\frac{1}{2} \|\Delta d\|_{W_2}^2,
\end{equation}
which is also a weighted-data residual as classic FWI. Similar with WRI, this means ESI may also suffer the local minima issue.
Eliminating $\lambda$ and $u$ gives the normal equation
\begin{equation}\label{eq:normalewi}
(S^H S + \beta B^H B) q = S^H d
\end{equation}
which is same as \citet[equation 10]{huang2018source}.
This equation shares the similarity with the normal equation \eqref{eq:normal} of WRI: the regularized Hessian matrix $S^H S + \beta B^H B$ has to be inverted for the solution of $q$, preventing the use of a time-domain wave equation solver. This prompted \citet{huang2018source} to conduct ESI in the frequency domain.

The above development puts classic FWI, WRI and ESI into the same framework using the Lagrangian-based adjoint state method. This allows us to better catch the similarity or differences among them. Table~\ref{tab:compare} lists the expressions of the misfit functional, the Lagrangian, the adjoint source and the normal equation of the three methods.

\begin{table}
  \caption{The misfit functional, Lagrangian, adjoint field and normal equations for FWI, WRI and ESI.}\label{tab:compare}
 \begin{tabular}{c|l|l|l|l}
  \hline
      & \textbf{Misfit} ($J$)            & \textbf{Lagrangian} ($\mathcal{L}$)         & \textbf{Adjoint field} $\lambda$                                    & \textbf{Normal equation} \\
  \hline
  FWI & $\frac{1}{2}\|Ru-d\|^2$                  & $J_\text{FWI}+\left\langle\lambda,Au-f\right\rangle$   & $S^{H}\Delta d$                                      & -  \\
  WRI & $\frac{1}{2}\|Ru-d\|^2+\frac{\beta}{2}\|q\|^2$  & $J_\text{WRI}+\left\langle\lambda,Au-f-q\right\rangle$ & $S^{H}\Delta d (=\beta q)$                & $(S^HS+\beta I)q=S^H(d-Sf)$ \\
  ESI & $\frac{1}{2}\|Ru-d\|^2+\frac{\beta}{2}\|Bq\|^2$ & $J_\text{ESI}+\left\langle\lambda,Au-q\right\rangle$   & $S^{H}\Delta d(=\beta B^HB q)$ & $(S^HS+\beta B^HB)q=S^Hd$   \\
  \hline
 \end{tabular}
\end{table}

\subsection{Role of the penalty parameter}\label{sec:beta}

Let us highlight the important role of the penalty parameter $\beta$ in WRI and ESI. The misfit functional is the summation over the data error $J^e$ and the penalty $J^p$ scaled by the penalty coefficient $\beta$:
\begin{equation}
  J_\text{WRI/ESI}=J^e + \beta J^p.
\end{equation}
\citet{fu2017discrepancy} have shown that
\begin{equation}
\frac{\mathrm{d} J^e}{\mathrm{d} \beta} \ge 0, \frac{\mathrm{d} J^p}{\mathrm{d} \beta} \le 0.
\end{equation}
Assuming a good initial guess, large $\beta$ will force the wave equation to be better satisfied, leading to smaller constribution to the misfit functional through $J^p$. Tuning $\beta$ should start from small values and increase through iterations to reach a rapid convergence \citep{fu2017discrepancy}. Based on the discrepancy principle, \citet{fu2017discrepancy} designed an automatic scheme to choose a proper $\beta$, starting the outer iterations for nonlinear optimization from $\beta=0$. It is noteworthy that the discrepancy principle applies equally well to WRI.

%
\section{Solving normal equation by conjugate gradient method}

As aforementioned, the presence of the regularized Hessian in the normal equation prevents a direct use of time stepping solvers for WRI and ESI.
\cite{wang2016full} provided the first workaround to conduct WRI into time domain. The idea is to use approximate wavefields by neglecting high order terms, so that the quantities are computable via the solution of wave equation using time stepping. Another attempt was made by \cite{rizzuti2020dual} based on a dual formulation. However, similar approximation has to be introduced in order to make the computation feasible, and twice of the modeling is required to compute all field quantities. Instead of solving the normal equation for the wavefield, \cite{gholami2020extended}  proposed to solve an equivalent normal equation for a data vector. Unfortunately, the key problem introduced by penalty formulation, i.e. the presence of the Hessian matrix in the normal equation, still stands. To proceed, they made drastic approximation by introducing a scalar to replace the inverse of the regularized Hessian, and assume the Green's function between two iterates are invariant. The augmented Lagrangian including a quadratic penalty can then be solved by the standard alternating direction method of multipliers (ADMM)  \citep{Boyd_2011_DOS}. These approaches are inaccurate from a mathematical point of view. To mitigate the nonlinearity and different sensitivities to the model parameter and the wavefield, the linear inversion inside a nonlinear inversion should be solved very accurately.

Since the Hessian matrix in the normal equation \eqref{eq:normalewi} is symmetric semi-positive definite (SPD), it is natural to consider it as a linear system
\begin{equation}\label{eq:axb}
\underline{A} \;\underline{x} = \underline{b},
\end{equation}
where $\underline{A}:= S^H S + \beta B^H B$.
We solve it by a linear conjugate gradient (CG) method in a matrix free fashion, as presented in Algorithm \ref{alg:CG}:
\begin{itemize}
\item The right hand side, $\underline{b}=S^Hd=A^{-H}R^Hd$, can be computed by first projecting the data $d$ back to the wavefield ($R^H d$), and then performing an adjoint modelling to get the adjoint field $b$.

\item The key to the CG method is how to perform the matrix vector product. Given an input vector $p$,  the matrix vector product,  $\underline{A}p = S^H S p + \beta B^H B p $, can be computed by summing over two parts. The first part, $S^H S p= A^{-H}R^H R A^{-1}p$, is obtained by performing two modellings: first do a forward modelling using $p$ as the source; then extract the data at receiver, and project it back to the wavefield; finally, use the extracted data to do an adjoint modelling. The evaluation of the second part, $\beta B^H B p$, is straightforward and could be added to the first part.

\item The solution of the unknown vector after the convergence of CG serves as the extended source: $\underline{x}:=q$.

\item In terms of the storage request, each CG iteration needs at least four vectors (the solution vector $\underline{x}$, direction vector $p$, residual vector $r$ and the matrix vector product $\underline{A}p$) to drive the algorithm into the next iteration. Here these vectors are built by the wavefields at every time step at every grid point, demanding huge amount of storage requirement. For the moment, we consider the strategy by disk IO, and discuss two possible alternatives after the numerical test.
\end{itemize}

The CG method does not require the explicit expression of $(B^H B)^{-1}$. The method converges for any $\beta\geq 0$. This is advantageous over some other possibilities such as Gauss-Seidel iteration and recursion assisted by surrogate function  \cite{aghamiry2020accurate}, see Appendix \ref{appendix}. Thanks to the matrix free formulation,  the above CG procedure can equally be applied to WRI by simply choosing $B=I$ with a different right hand side. 
In terms of equation \eqref{eq:normalwri}, the right-hand-side term of the normal equation for WRI is $b=S^H(d-Sf)=A^{-H}R^{H}(d-RA^{-1}f)$, which can be computed by first extracting the forward wavefield excited by the true source $f$, then use the data difference ($d-RA^{-1}f$) as the virtual source to compute the adjoint wavefield $b$. Both time and frequency domain wave equation solvers can be applied when solving normal equation using CG.

\begin{algorithm}[H]
  \caption{Conjugate gradient algorithm for solving $\underline{A}\;\underline{x}=\underline{b}$ ($\underline{A}^H=\underline{A}$)}\label{alg:CG}
  \begin{algorithmic}[1]
  \State $x_0:=0$
  \State $r_0:=b-\underline{A} x_0=b$
  \State $p_0:=r_0$
    \For{$k=0,\cdots,N_{CG}-1$}
    \State $\alpha_k=\frac{r_k^H r_k}{p_k^H \underline{A}p_k}$
    \State $x_{k+1}=x_k + \alpha_k p_k$
    \State $r_{k+1}=r_k - \alpha_k \underline{A} p_k$
    \If{$\|r_{k+1}\|< tol\cdot \|r_0\|$ or $k==N_{CG}-1$}
    \State      exit loop\;
    \EndIf
    \State $\beta_{k} = \frac{r_{k+1}^H r_{k+1}}{r_k^H r_k}$\;
    \State $p_{k+1}=r_{k+1} + \beta_k p_k$\;
    \EndFor
  \end{algorithmic}
\end{algorithm}

A special case comes with $\beta=0$. The normal equation then becomes
\begin{equation}
  S^H S q = S^H d .
\end{equation}
When $q$ is restricted to the true source location, finding $q$ based on the given data $d$ is exactly the source wavelet estimation process, which is solvable using CG since $S^H S$ is SPD. In this case, the problem degenerates to classic FWI: we simply use the observed data $d$ to estimate a wavelet $q$ at the true source location, then minimize the misfit functional iteratively.

The optimization is consequently performed using the variable projection method \citep{golub1973differentiation}
due to the dramatic difference in the sensitivity of the misfit functional to the extended sources $q$ and the model parameters $m$.
There are nested loops with the nonlinear outer loop for searching model parameters $m$ and the linear inner loop for
updating extended source $q$ \citep{huang2018source}. In the frequency formulation, \citet{aghamiry2019improving} has proposed
a linear outer loop over the wavefield $u$ and a linear inner loop over specific model parameters $m$. He has taken
benefit from the linear property of the scalar wave equation with respect to the square of slowness for fixed wavefield $u$.
Still the inner loop has to be solved through direct solver or iterative solver in this frequency formulation. In our strategy
for time formulation, as sketched in Algorithm \ref{alg:nestedloop}, the outer loop is the iterative update for model parameters $m$,
while the inner loop is the solution of the normal equation with respect to the extend source $q$. In order to mitigate the potential
non-linearity of the inverse problem, the inner loop for the linear optimization, {\it i.e.} the CG iteration, should converge
with good precision. The recursion strategy illustrates that the inner loop should be solved quite accurately \citep{aghamiry2019improving}.

\begin{algorithm}[!tbh]
  \caption{Nested loops for variable projection between $m$ and $q$}\label{alg:nestedloop}
  \begin{algorithmic}[1]
  \State Given initial model: $m_0$\;
  \For{$i=0,\cdots,N_{outer}-1$ (outer loop: $\min_m J_\text{ESI}(m, q)$)}
  \State solve $ (S^H S + \beta B^H B) q= S^H d$ using Algorithm \ref{alg:CG}
  (inner loop: $\min_q J_\text{ESI}(m, q)$)
    \State compute forward field: $u=A^{-1}q$\;
    \State compute adjoint field: $\lambda=\beta B^H B q$\;
    \State calculate gradient: $\partial_m J=\sum_\omega \lambda^H \frac{\partial A}{\partial m} u $\;
    \State estimate descent direction: $\Delta m=- \tilde{H}^{-1} \partial_m J$ ($\tilde{H}$=approximate Hessian by L-BFGS algorithm)\;
    \State find a good step size $\alpha$ by line search to update model\;
    $m_{i+1}:=m_{i}+\alpha \Delta m$\;
  \EndFor
  \end{algorithmic}
\end{algorithm}

\section{Numerical tests}

We now numerically access the feasibility of time-domain WRI and ESI using CG method. Schematically, we always first invert for the extended source, then deduce the forward and adjoint fields to build the gradient of the misfit functional. It should be highlighted that the first stage in conventional WRI actually focuses on inverting the wavefield $u$ instead of extended source $q$, which is different from what the WRI computing scheme presented here. In all the numerical experiments as follows, we use 10 CG iterations to solve the normal equation for WRI and ESI. The annihilator is chosen to be $B=\|x-x_s\|$. Note that we add a small positive number to $B$ to ensure $B$ is not singular around the source $x_s$.

\subsection{Data fitting in reflection case}

\begin{figure}[!tbh]
 \centering
 \includegraphics[width=\linewidth]{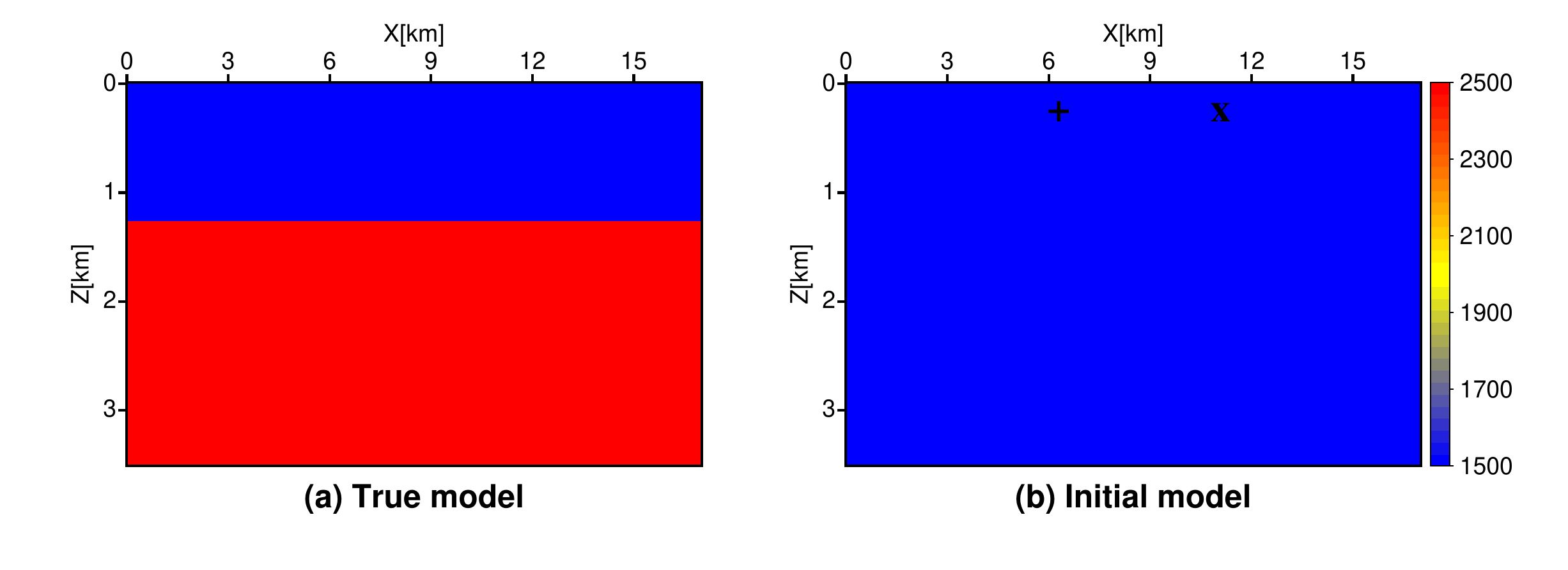}
 \caption{(a) True velocity model including two layers of velocities 1500 m/s and 1800 m/s; (b) Initial velocity model of constant velocity 1500 m/s, where the source and receiver positions are marked by  $+$ and $\times$.}\label{fig:twolayer}
\end{figure}

We start with a 2D synthetic model including two layers of different velocities (1500 m/s and 1800 m/s), as shown in Figure~\ref{fig:twolayer}a. The initial model displayed in Figure~\ref{fig:twolayer}b is homogeneous with velocity equal to the first layer of the true model. 
The locations of the source and the receiver has been marked up in Figure~\ref{fig:twolayer}b. The purpose of this test is to understand the data fitting and the gradient of FWI, WRI and ESI. 
The modelling scheme after finite difference discretization is $\mathcal{O}(\Delta t^2, \Delta x^8)$. An absorbing boundary condition is used to get rid of the free surface and boundary reflections.

\begin{figure}[!tbh]
 \centering
 \includegraphics[width=\linewidth]{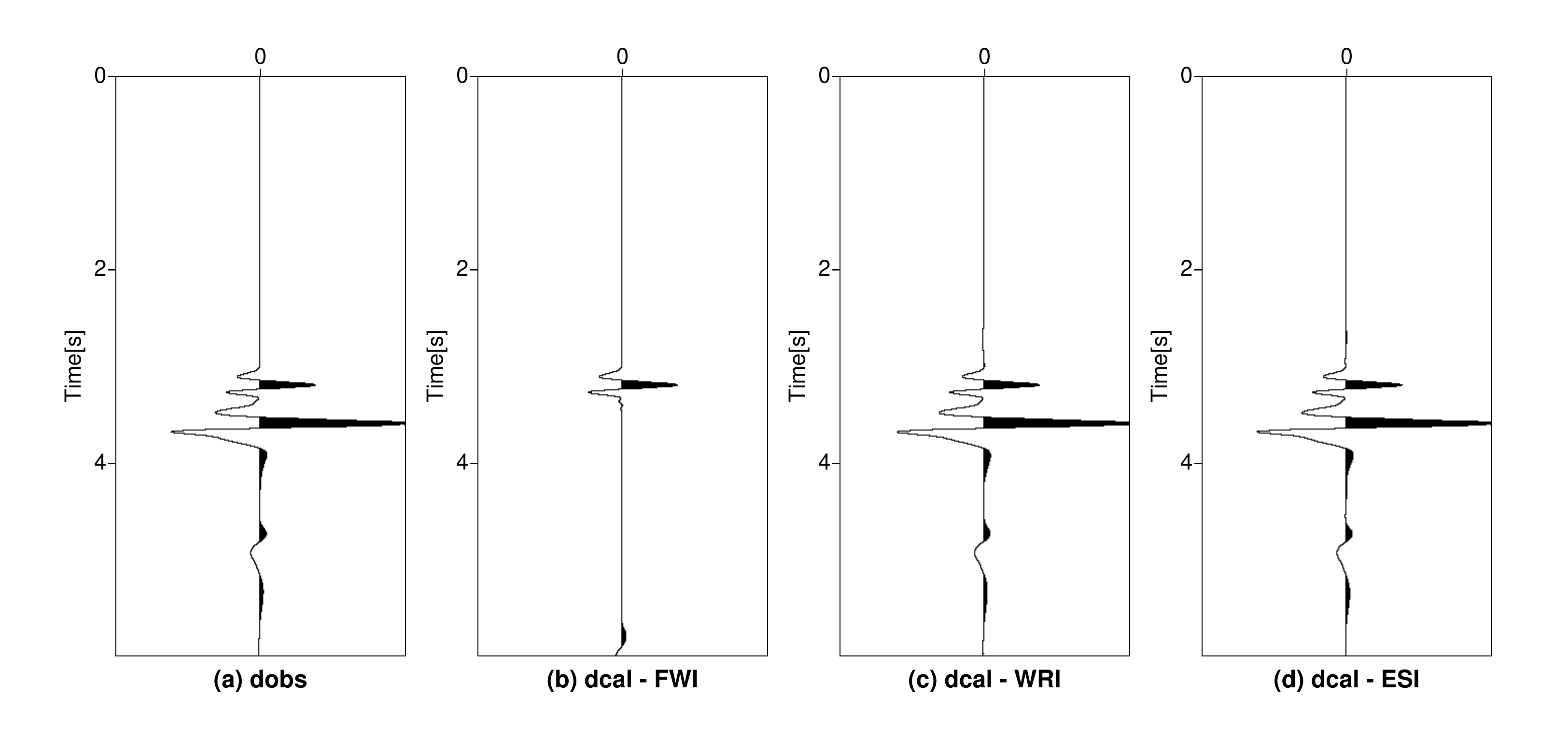}
 \caption{(a) Observed data simulated from true model; (b) Synthetic data by classic FWI (direct arrival from source to receiver is reproduced, but reflection from the interface is missing); (c) Synthetic data by WRI; (d) Synthetic data by ESI.
}\label{fig:dat_compare}
\end{figure}

The observed data from the true model has been displayed in Figure~\ref{fig:dat_compare}a. Since the velocity of the initial model is the same as the upper layer of the true model, the synthetic data in classic least-squares FWI reproduces direct arrivals of the observed data perfectly; however, the reflections from the interface between two layers are completely missing, as can be seen from Figure~\ref{fig:dat_compare}b. Figures~\ref{fig:dat_compare}c and \ref{fig:dat_compare}d shows that both WRI and ESI generate data very close to observed data, mimicking that the simulation is performed in the true model such that the data can be fitted almost perfectly. By domain extension and recursive data assimilation, the information from the observed data are added back to the modelled data, behaving that the true physics has been taken into account. As displayed in Figure~\ref{fig:dat_compare}a, classic FWI only considers information in the starting model, thus is blind to the information in the data. 

\begin{figure}[!tbh]
 \centering
 \includegraphics[width=\linewidth]{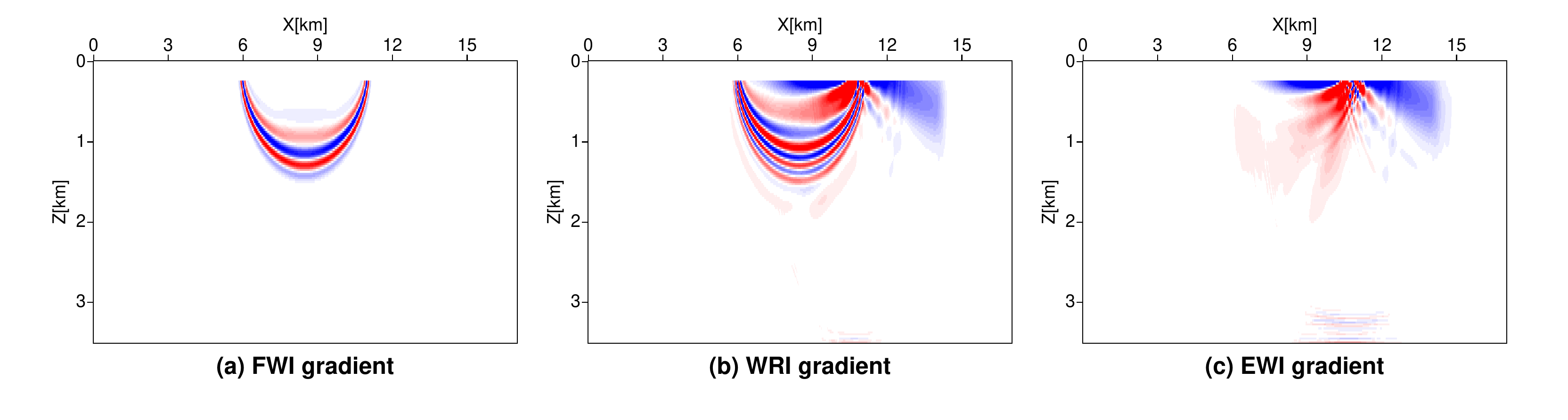}
 \caption{(a) Classic FWI, (b) WRI and (c) ESI gradients representing the sensitivity of each method to model variations. ESI gradient presents one-sided illumination, significantly different from the two-sided FWI and WRI gradients.}\label{fig:grad}
\end{figure}

The gradients from FWI, WRI and ESI have been plotted in Figure~\ref{fig:grad}. The FWI gradient in Figure~\ref{fig:grad}a shows better source receiver symmetry in the banana donut of the energy isochrone, while the ESI gradient in Figure~\ref{fig:grad}c shows more energy comes from receiver side, including a lot of low frequency energy in the first Fresnel zone. Figure~\ref{fig:grad}b is in between the previous two, with more energy on the receiver side than on the source side.  This is quite understandable as WRI uses both the physical source and the extended source. In both WRI and ESI, the extended source extracts information based on the data from the receiver location.

\subsection{Misfit landscape in transmission case}

\begin{figure}[!htb]
 \centering
 \includegraphics[width=0.75\linewidth]{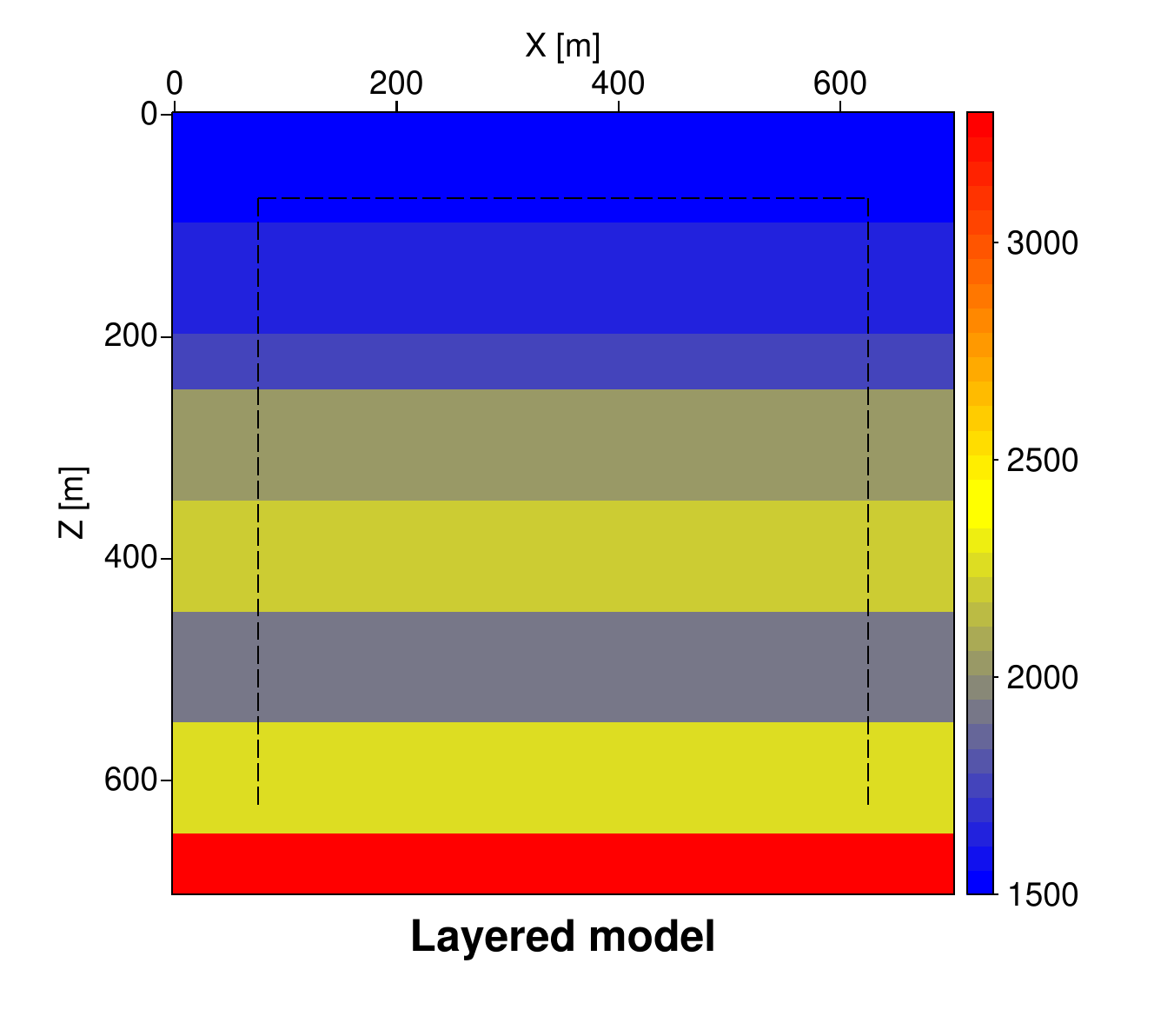}
 \caption{The 1D layered velocity model used for generating the initial model to study the variation of the misfit functional with respect to velocity errors. The dashed line indicates the location of the receivers, while the source coincides with the receiver in the top left corner.}\label{fig:vp_true_scan}
\end{figure}

\begin{figure}[!htb]
 \centering
 \includegraphics[width=0.75\linewidth]{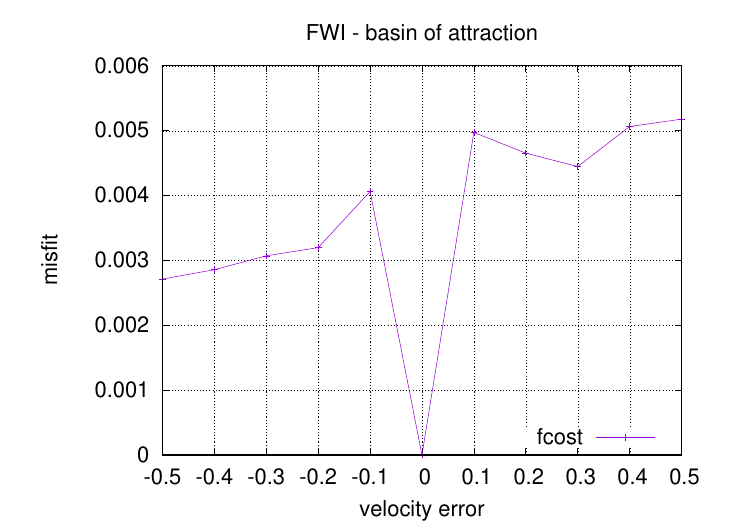}
 \caption{The variation of the misfit functional with respect to velocity errors for classic FWI}\label{fig:misfit_vp_fwi}
\end{figure}

\begin{figure}[!htb]
 \centering
 \includegraphics[width=0.75\linewidth]{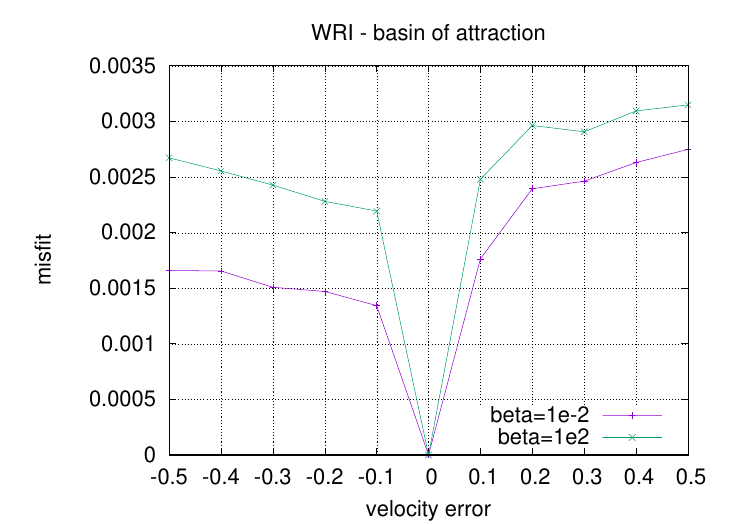}
 \caption{The variation of the misfit functional with respect to velocity errors for WRI}\label{fig:misfit_vp_wri}
\end{figure}

\begin{figure}[!htb]
 \centering
 \includegraphics[width=0.75\linewidth]{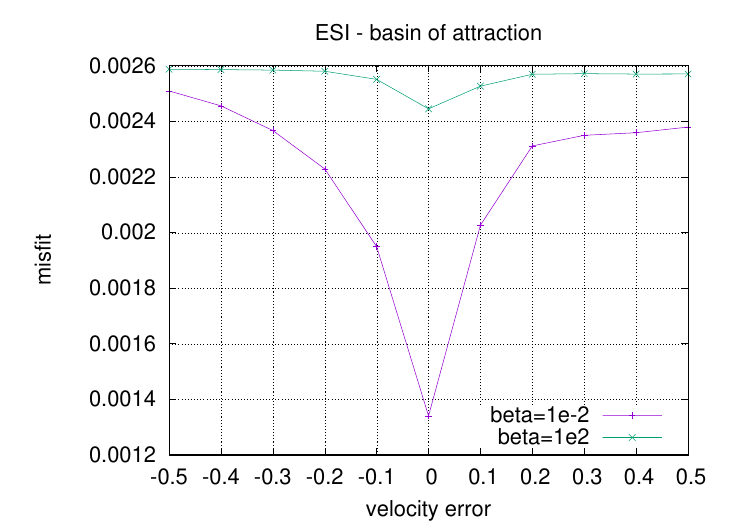}
 \caption{The variation of the misfit functional with respect to velocity errors for ESI}\label{fig:misfit_vp_ewi}
\end{figure}

In order to better understand  the capability of WRI and ESI methods in getting rid of local minima inherent in classic FWI, we now study the landscape of the misfit functional, based on a simple 1D layered model shown in Figure~\ref{fig:vp_true_scan}. We consider the surface acquisition plus two vertical lines of receivers distributed on two sides of the model. The source has been deployed at the top left corner coinciding with the receiver position, such that there are sufficient transmission energy.
We first generate a number initial velocity model according to the formula $m_0=(1+\epsilon)m_{true}$ ($-0.5\leq \epsilon \leq 0.5$), where $m_{true}$ is the true layered model.
We then compute the misfit functions by taking different values of the velocity error $\epsilon$. The basin of the attraction of FWI, WRI and ESI are shown in Figures~\ref{fig:misfit_vp_fwi}, \ref{fig:misfit_vp_wri} and \ref{fig:misfit_vp_ewi}. We examine the WRI and ESI using a large penalty $\beta=100$ and a small penalty $\beta=0.01$. The classic FWI misfit plotted in Figure~\ref{fig:misfit_vp_fwi} shows multiple local minimum, while the convexity domain is rather narrow. The  WRI misfit in Figure~\ref{fig:misfit_vp_wri} has larger convexity domain than FWI, but still suffers from local minimum issue, especially when $\epsilon>+0.2$ for $\beta=100$. Figure~\ref{fig:misfit_vp_ewi} confirms that ESI further extends the search space, possessing wider basin of attraction, compared to WRI and least-squares FWI. It is interesting to note that due to the absence of the knowledge on the physical source, the observed data can never be fit perfectly. Hence, the ESI misfit is not zero even if the input model has no error compared to the true model. For both WRI and ESI, the landscape of misfit functional becomes flattened when the penalty is small. This indicates that the penalty must be judiciously chosen, increasing from a small value with wider basin of attraction of the misfit, to a large value for fast convergence when the model gets close to the ground truth.

\subsection{Marmousi-II test: preliminary validation of time-domain ESI}

\begin{figure}[!htb]
 \centering
 \includegraphics[width=\linewidth]{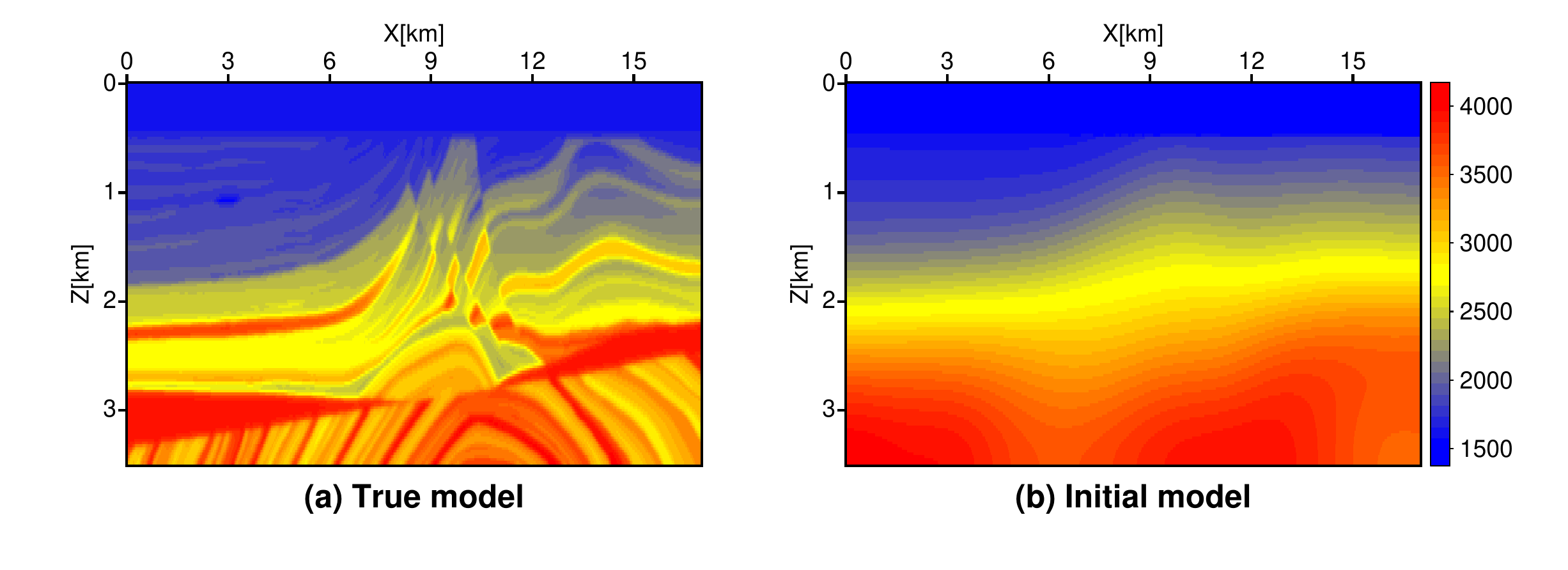}
 \caption{Marmousi-II test: (a) True model; (b) initial model. }\label{fig:marmousi}
\end{figure}

Finally, we validate our time-domain ESI implementation based on the Marmousi-II model. As shown in Figure~\ref{fig:marmousi}b, the starting model has been highly smoothed based on the  true model in Figure~\ref{fig:marmousi}a. The model has the water depth of 500 m. We consider a surface acquisition geometry: 28 sources are evenly distributed at 25 m of the water depth, while 641 receivers are uniformly deployed at the same depth. The observed data was generated using a Ricker wavelet of 5 Hz peak frequency. The tuning of the penalty parameter $\beta$ was automatically determined through iterations using discrepancy principle.

\begin{figure}[!htb]
 \centering
 \includegraphics[width=0.7\linewidth]{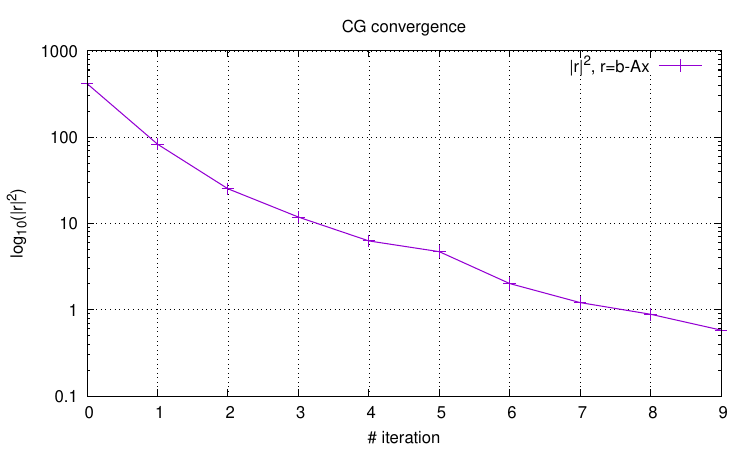}
 \caption{The convergence history of CG for ESI in the 1st CG loop.}\label{fig:cg_conv}
\end{figure}

\begin{figure}[!htb]
 \centering
 \includegraphics[width=\linewidth]{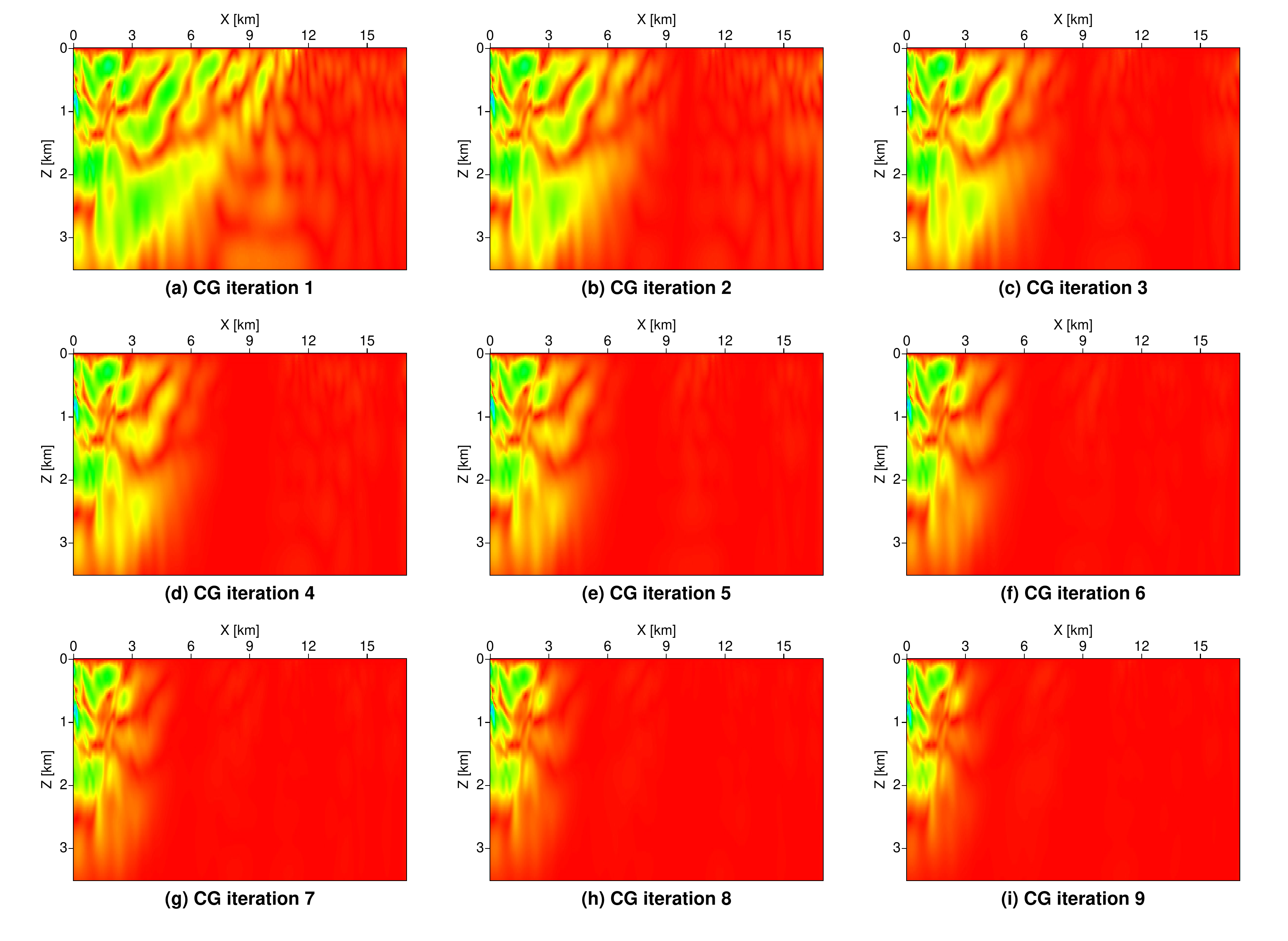}
 \caption{The evolution of extended source at 2 Hz for the 1st shot in the 1st CG loop using discrete Fourier transform. Red: amplitude=0; Green: high amplitude.}\label{fig:extendedsrc}
\end{figure}

We observed that our CG algorithm enjoys a fast convergence rate, as can be seen from Figure~\ref{fig:cg_conv}. Figure~\ref{fig:extendedsrc} computed the amplitude of the extended source at 2 Hz for the 1st shot in the 1st CG loop,  using discrete time Fourier transform (see next section).
It illustrates that the extended source varies dramatically at early iterations, while reach a steady state after a number of CG iterations: the energy of extended source indeed focuses quickly to the true source location (the 1st shot resides in the top left part of the model) through CG iterations, even though the synthetic data in ESI is purely deciphered from the recursion process according to the observed data without the knowledge of the true physical source.

\begin{figure}[!htb]
 \centering
 \includegraphics[width=\linewidth]{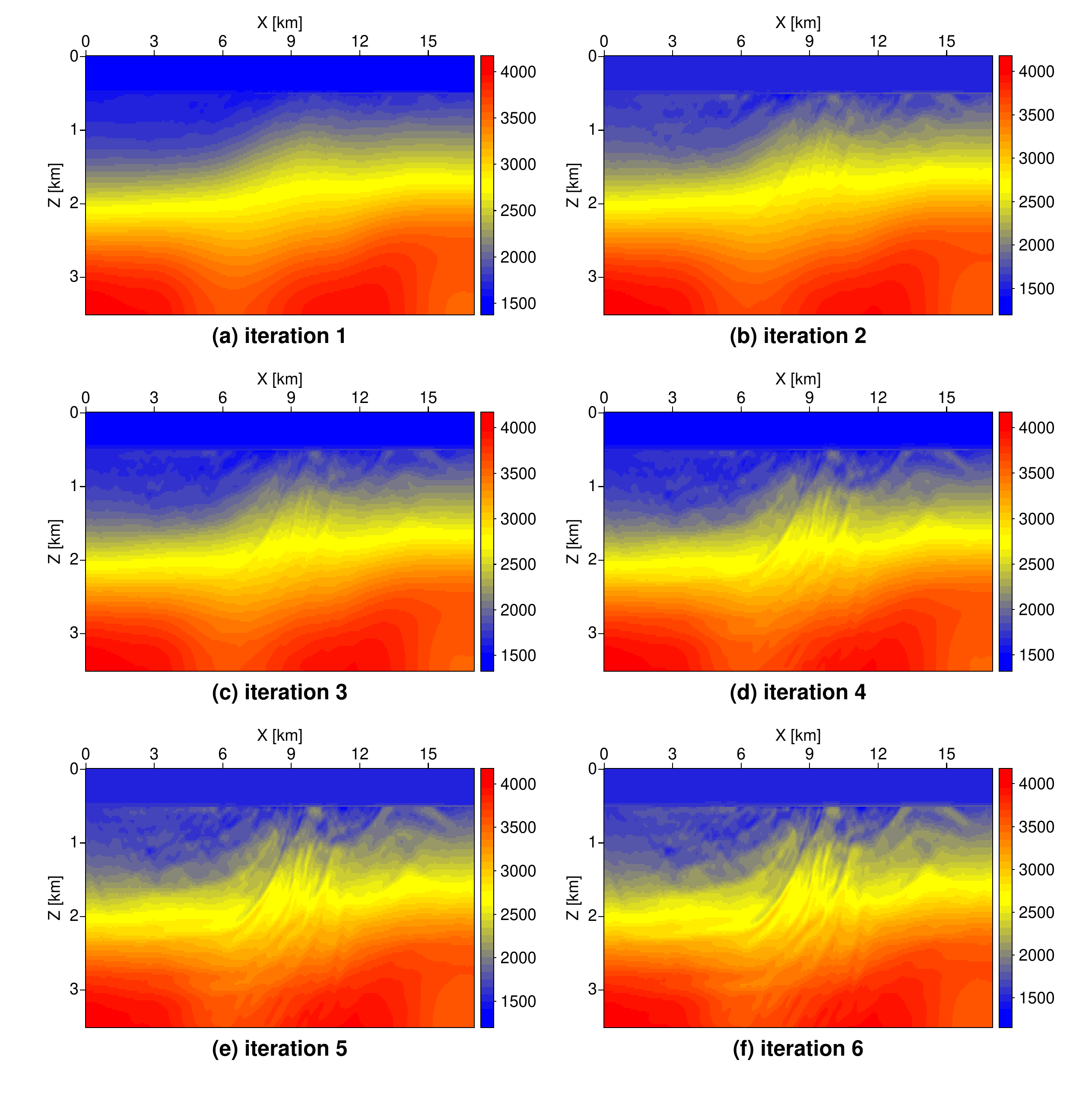}
 \caption{The model updated at different iterations. }\label{fig:iterates}
\end{figure}

We managed to complete 6 ESI iterations. The updated models after limited number of iterations are displayed in Figure~\ref{fig:iterates}, showing that the model is evolving to the true model. Unfortunately, the whole inversion still converges very slowly so that we have to terminate it on the way, since the 6 updates of the 2D test has consumed more than 72 hours of CPU time due to intensive computational cost and slow IO in order to read and write huge amount of wavefield snapshot to the disk. This motivates us to focus on the development of low storage ESI approach without disk IO.

\section{Mitigating the storage demand}

Unlike frequency-domain implementations, WRI and ESI in the time domain pose significant challenges on the memory requirement, because all snapshots of the wavefield simulations have to be stored (RAM or disk), before the two parts of the matrix vector product are assembled. Consequently, this huge storage requirement becomes one of the major bottlenecks to push forward time-domain WRI and ESI to practical applications. Other bottlenecks such as nested loop for expensive computation and the need for good preconditioning are out of scope of the discussion for the moment. A natural idea is to compress the wavefield, which is still not practical for 3D or large-scale applications: Most of the lossless compression technique cannot achieve an extremely high compression ratio such that the wavefield snapshots can be saved in the RAM. Our numerical experience also shows that the lossy compression technique breaks the orthogonality of the residual vectors between two iterates, and thus CG algorithms do not work. We remind that switching to Gauss-Seidel iterations and the recursion scheme assisted by surrogate function (cf. Appendix \ref{appendix}) will encounter similar memory issues.

While the time-stepping modeling is stilled, we have made two attempts to bypass the storage issue: the data space approach and the frequency-domain inversion by time domain modelling. Although both solutions are not very successfull, we discuss them here to highlight the significant challenges in handling the memory issue, hoping to draw the attention of the community to the memory issues of extended domain inversion methods.

\subsection{Switching from the wavefield to data space}

Let us define the synthetic data extracted at receiver location as $w:=Ru$. Then we have
\begin{equation}
  \begin{split}
  w=& Ru=RA^{-1}q \quad (u=A^{-1}q \mbox{ from equation \eqref{eq:ewifwd}})\\
  =&Sq=S(\beta B^H B)^{-1} \lambda \quad (q=(\beta B^H B)^{-1}\lambda \mbox{ from equation \eqref{eq:ewilq}})\\
  =&S(\beta B^H B)^{-1} A^{-H} R^H (d-Ru) \quad (\lambda=A^{-H} R^H(d-Ru) \mbox{ from equation \eqref{eq:ewiadj}})\\
  =&S(\beta B^H B)^{-1} S^H (d-w) \quad (S=RA^{-1}, w=Ru)
  \end{split}
\end{equation}
yielding
\begin{equation}\label{eq:wd}
    (I +S(\beta B^H B)^{-1} S^H ) w = S(\beta B^H B)^{-1} S^H d
\end{equation}
Define
\begin{equation}\label{eq:vd}
  v:=S(B^H B)^{-1} S^H d.
\end{equation}
Equation \eqref{eq:wd} then becomes
\begin{equation}\label{eq:wv}
(\beta I +S(B^H B)^{-1} S^H )w = v.
\end{equation}

The data vector $v$ can be computed by first performing an adjoint modelling using $d$ as the adjoint source, then apply the operator $(B^H B)^{-1}$ on the adjoint field, and finally extract the data at receiver locations after a forward modelling using scaled adjoint field. Thanks to the symmetric positive definite matrix $(\beta I +S(B^H B)^{-1} S^H )$, the solution of linear system for $w$ in \eqref{eq:wv} can also be found using conjugate gradient method in a matrix free manner. Once the synthetic data $w$ is computed, the adjoint field can be computed according to equation \eqref{eq:ewiadj}. Then the extended source is obtained from equation \eqref{eq:ewilq}: $q=(\beta B^H B)^{-1}\lambda$. The field $u$ is immediately available after a forward modelling based on the extended source $q$: $u=A^{-1} q$. Finally, the ESI gradient is built by cross-correlation between the forward and the adjoint fields via  \eqref{eq:gradewi}.

In fact, the above computing workflow can also be discovered according to Sherman-Morrison-Woodberry (SMW) formula \citep[p.51, equation 2.1.4]{Golub_1996_MATCOMP}
\begin{equation}
  (P+UV^H)^{-1}= P^{-1} - P^{-1}U(I_k+V^H P^{-1}U)^{-1}V^H P^{-1},
\end{equation}
where $P\in \mathbb{R}^{n\times n}$ is invertible while $U\in \mathbb{R}^{n\times k}$, $V\in \mathbb{R}^{n\times k}$.
Assume $\beta B^H B$ is invertible. We can express the inverse of the normal operator in the following:
\begin{equation}
   (S^H S + \beta B^H B)^{-1} 
  =(\beta B^H B)^{-1} - (\beta B^H B)^{-1} S^H (I + S(\beta B^H B)^{-1}S^H)^{-1} S (\beta B^H B)^{-1}.
\end{equation}
The extended source is then
\begin{equation}
  \begin{split}
    q=&(S^H S + \beta B^H B)^{-1} S^H d\\
    =&(\beta B^H B)^{-1}  S^H d- (\beta B^H B)^{-1} S^H (I + S(\beta B^H B)^{-1}S^H)^{-1} S (\beta B^H B)^{-1} S^H d\\
    =&(\beta B^H B)^{-1}  S^H \Big(d- \underbrace{ (\beta I + S(B^H B)^{-1}S^H)^{-1} \overbrace{ S ( B^H B)^{-1} S^H d}^v}_w\Big)
  \end{split}
\end{equation}
 The dimensionality of $w$ living in data space is much lower than the dimensionality of the wavefield space in the time domain. It can therefore be directly stored in memory. Finally, the extended source will be given by
\begin{equation}\label{eq:q}
  q = (\beta B^H B)^{-1} S^H (d-w).
\end{equation}
Since the synthetic data $w$ reads
\begin{equation}
  w:= (\beta I + S(B^H B)^{-1}S^H)^{-1} S( B^H B)^{-1}  S^H d,
\end{equation}
it is therefore natural to see that the synthetic data in ESI is able to fit the observed data precisely as $w\rightarrow d$ ($\beta\rightarrow 0$).

A pure time domain implementation of the above algorithm is still difficult due to the need to access the adjoint field as the source for forward modelling (applying operator $S$ after application of $(B^H B)^{-1}S^H$).
The scaled adjoint field has opposite direction in time during the simulation compared to the forward field. This creates also challenges due to gigantic memory overhead and immense amount of IO traffic, particularly for 3D practical applications. In fact, building the FWI gradient faces up to the same issue as the recursion process. For ESI in non-attenuating medium, we can reconstruct the  wavefield on the fly by decimation and interpolation of the stored boundaries \citep{Yang_2016_WRB}. For ESI in attenuative medium we may have to resort to the optimal checkpointing \citep{Symes_2007_RTM} or CARFS algorithm  \citep{Yang_2016_CAR}, since recomputing wavefield by decimation and interpolation of the stored boundaries is no more stable.

The above computational strategy applies equally to WRI: what we need is just replacing $d$ with $d-Sf$ and set $B=I$ in computing $v$ and $q$.  That is,
\begin{equation}
  \begin{split}
    q=&(S^H S + \beta I)^{-1} S^H (d-Sf)\\
    =&\beta^{-1}  S^H (d-Sf)- \beta^{-1} S^H (I + \beta^{-1}SS^H)^{-1} \beta^{-1} SS^H (d-Sf)\\
    =&\beta^{-1}  S^H \Big(d-Sf- \underbrace{ (\beta I + SS^H)^{-1} \overbrace{ S S^H (d-Sf)}^v}_w\Big)
  \end{split}
\end{equation}
where
\begin{equation}\label{eq:wriw}
 w =(\beta I + S S^H)^{-1} S S^H(d-Sf),
\end{equation}
yielding
\begin{equation}
  \begin{split}
  w = &\beta^{-1}S(B^H B)^{-1} S^H(d-Sf-w)\\
  =&S(\beta B^H B)^{-1} S^H(d-Sf-w)\\
  =&Sq =RA^{-1}(Au-f) \\
  =&Ru -Sf
  \end{split}
\end{equation}
so that the synthetic data in WRI is $Ru=w+ Sf$, and the data residual is $\Delta d= d-Ru = d-Sf -w$ which can also be very small as $w\rightarrow (d-Sf)$ if $\beta\rightarrow 0$ according to \eqref{eq:wriw}. The wavefield $u$ in WRI is obtained from the forward modelling based on the source $f+q$ instead of only $q$.

\subsection{Frequency domain inversion by time domain modelling}

The second idea is to do  frequency domain inversion with time domain modelling, inspired by  \citet{Sirgue_2008_FDW}. The main merit of this approach is that all (forward and adjoint) wavefields needed are stored at few frequencies, while the time domain wavefield can be computed on the fly during time stepping, via discrete time Fourier transform (DTFT):
\begin{equation}\label{eq:dtft}
u(x,f_m) = \sum_{n=0}^{N_t-1} u(x,n\Delta t) \exp(j 2\pi f_m n\Delta t), \quad m=1,\cdots,N_f,
\end{equation}
where $f_m$ refers to the $m$-th discrete frequencies.
The modelling of the adjoint field in the time domain requires the adjoint source  (the data residual) as a time series, meaning that all the frequencies are needed. It is unfortunate that the data and the wavefield are only accessible at limited number of discrete frequencies. To supply such an adjoint source $s$, we use a time series computed by inverse DTFT based on few discrete frequencies available:
\begin{equation}\label{eq:idtft}
s(x,n\Delta t) =\frac{1}{N_t} \sum_{m=1}^{N_f} s(x,f_m) \exp(-j 2\pi f_m n\Delta t), \quad n=0,\cdots,N_t-1.
\end{equation}
Note that the number of discrete frequencies $N_f$ is much less than total number of time steps $N_t$ in the simulation. The stringent storage demand is therefore resolved.

In fact, equation \eqref{eq:idtft} should be understood as the adjoint of the DTFT scaled by the factor $1/N_t$, which makes the CG feasible. One may consider apply  the same scaling factor $1/\sqrt{N_t}$ for both forward DTFT and its inverse, to make it fully adjoint. This is in fact unnecessary as each CG iteration requires one matrix vector product. The matrix vector product operates on a frequency domain wavefield, which is first transformed into time domain by inverse DTFT using \eqref{eq:idtft}, then performs time domain adjoint and forward modelling, and finally the  forward field are transformed back into frequency domain using DTFT \eqref{eq:dtft}. Although the constant scaling factor $1/N_t$ only appears in inverse DTFT, the resulting matrix vector product will be the same as normalizing DTFT and inverse DTFT with the same factor $1/\sqrt{N_t}$.

Let us consider a frequency domain formulation of ESI.
\begin{equation}
  J = \frac{1}{2}\|\bar{d}-F Ru\|^2 + \frac{\beta}{2}\|Bq\|^2
  = \frac{1}{2}\|\bar{d}-F S F^H \bar{q}\|^2 + \frac{\beta}{2}\|BF^H \bar{q}\|^2 
\end{equation}
such that the normal equation is given by
\begin{equation}
( (FS F^H )^H FSF^H + \beta (BF^H)^H BF^H) \bar{q} = (FSF^H)^H \bar{d}
\end{equation}

Due to the fact that $F F^H = I^{N_f\times N_f}$ and the spatial operator $B$ and the temporal operator $F$ commute, we immediately reduce the above normal equation to
\begin{equation}
( F S^H F^H FSF^H + \beta B^H B) \bar{q} = FS^H F^H \bar{d}
\end{equation}

Let us denote DTFT in \eqref{eq:dtft} as $F\in\mathbb{C}^{N_f\times N_t}$ and the inverse DTFT in \eqref{eq:idtft} as $F^H\in \mathbb{C}^{N_t\times N_f}$.
The solution of normal equation using CG can then be understood as applying an additional DTFT operator on both sides of the equation \eqref{eq:normalewi}:
\begin{equation}\label{eq:normaldtft}
F(S^H S +\beta B^H B) F^H \bar{q} = FS^H F^H \bar{d}
\end{equation}
where the time domain wavefield and the observed data in \eqref{eq:normalewi} have been represented by frequency domain counterpart through inverse DTFT, that is,
\begin{displaymath}
  q(x,t) = F^H \bar{q}(x,f), \quad d(x_r, t)=F^H \bar{d}(x_r, f).
\end{displaymath}
It is clear that $F(S^H S +\beta B^H B) F^H $ is SPD, and therefore can be solved using CG. In particular, when the number of frequencies are the same as the number of time steps  ($N_f=N_t$), the solution of normal equation in \eqref{eq:normaldtft} using CG becomes standard preconditioned CG \citep[section 5.1 preconditioning]{Nocedal_2006_NOO}.

Once the extended source is obtained, we immediately obtain the gradient of the misfit functional:
\begin{equation}
  \partial_m J_\text{ESI} =\Re \sum_\omega \lambda^H \frac{\partial A}{\partial m} u,
\end{equation}
where the forward field $u=A^{-1}q$ should be computed by a forward modelling using the estimated extended source after CG, while the adjoint field is given by $\lambda=\beta B^H B q$. A pure frequency domain implementation through the solution of Helmholtz equation has been done by \citet{huang2018source}, while frequency domain inversion using time domain modeling is another alternative to achieve the same goal thanks to DTFT \citep{Sirgue_2008_FDW}. 

\begin{figure}[!htb]
 \centering
 \includegraphics[width=0.85\linewidth]{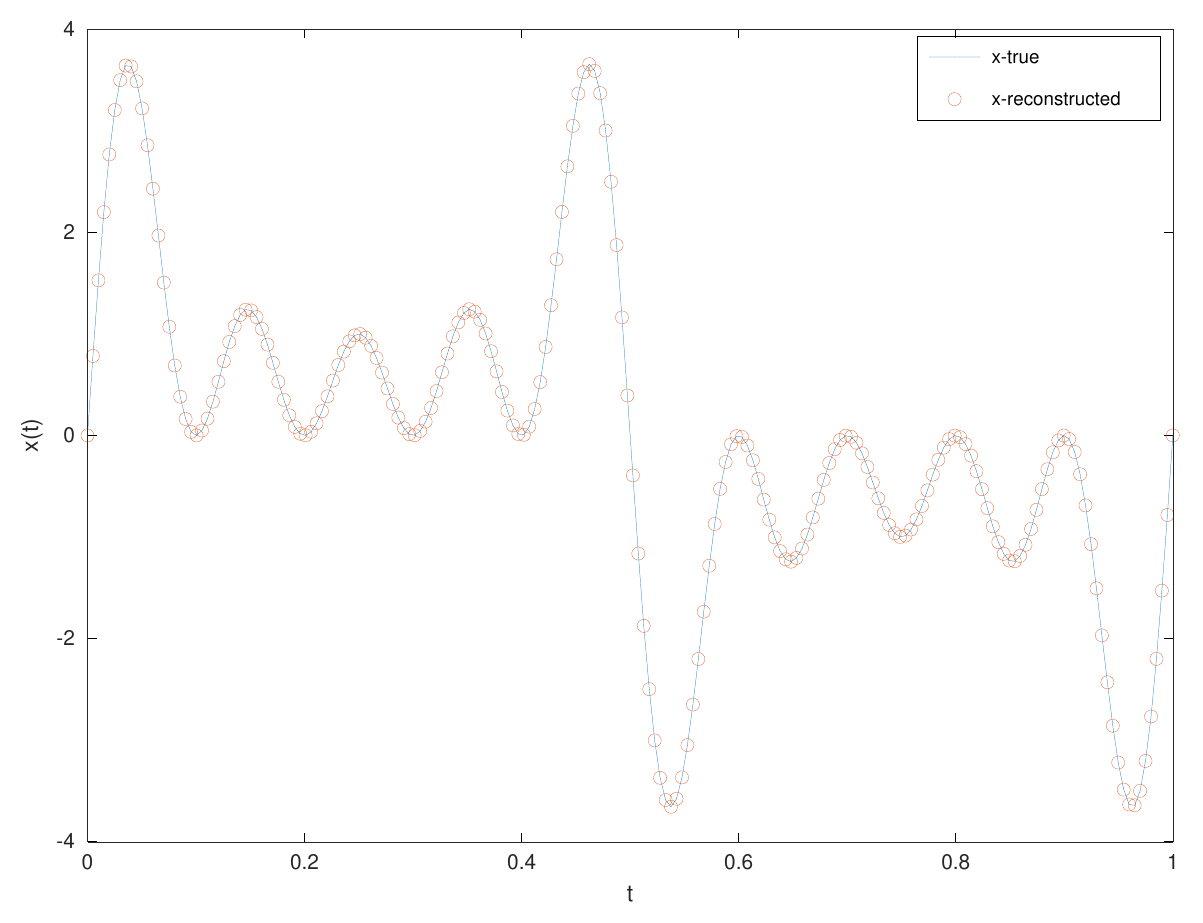}
 \caption{Perfect reconstruction of $x(t)$ of temporal length $N_t=200$ using DTFT and CG of length $N_f=5$ in frequency domain.}\label{fig:cg_dtft}
\end{figure}

As a proof of concept, we generate a time domain signal $x(t)=\sum_{m=1}^{N_f} \sin(2\pi f_m t), t\in [0,1]$ up to $N_t=200$ temporal samples, where $N_f=5$ frequencies are specified, $[f_1,f_2,f_3,f_4,f_5]=[1,3,5,7,9]$ Hz. A Toeplitz matrix of size $N_t\times N_t$ has been created to mimic the role of regularized Hessian operator in the time domain. We can then construct the right hand side of linear system $b=\underline{A}x$. By converting $b$ into frequency domain using DTFT, we use CG to solve the linear inverse problem to obtain a frequency domain signal which is only of length $N_f=5$ rather than $N_t=200$. Finally, we convert the frequency domain solution back to time domain. The reconstructed signal matches the true one perfectly in the time domain, as can be seen in Figure~\ref{fig:cg_dtft}.

\section{Conclusion}

This paper formulates classical FWI, WRI and ESI in a unified framework using the Lagrangian-based adjoint state method.
We show both WRI and ESI can be expressed as a weighted norm of least-squares FWI, hence the two methods potentially suffer the local minima issue.
Formulating the wave equation as a penalty term in the misfit functional leads to the presence of the regularized Hessian in the normal equations of WRI and ESI, preventing a direct use of time stepping modelling schemes  for wavefield simulations. We have considered the conjugate gradient method to efficiently solve the normal equation, allowing both time and frequency domain wave equations to be used in a matrix free fashion. The numerical tests show that ESI has larger convexity domain compared with WRI and classic FWI. A pure time-domain implementation is feasible from the algorithmic point of view, but faces the memory bottleneck. We made two additional attempts to alleviate this issue which requires further investigations.
Therefore, to make WRI or ESI a viable candidate for practical application, smart strategies must be designed to mitigate the storage requirement for time-domain inversion, and efficient preconditioner should be developed to reduce the number of inner CG iterations in the nested loop.

\section*{Acknowledgments}
We thank the inspiring communications with Jean Virieux and William W. Symes, in particular on the understanding of WRI and ESI approaches.

\appendix

\section{Another two methods for time-domain ESI}\label{appendix}

\subsection{Gauss-Seidel iterations}

The relations in equations \eqref{eq:ewiadj}, \eqref{eq:ewilq} and \eqref{eq:newq2} derived using adjoint state method give us several coupled equations equivalent to the normal equation \eqref{eq:normalewi}, but shows possibility to compute wavefield by time stepping wave equation solvers if the right hand side is known.
It inspires us to consider a recursive strategy to solve the normal equation, from a starting guess of the wavefield solution $u^0$. We can then recursively compute a new iterate of the field $\lambda^{k}$ by \eqref{eq:ewiadj}, hence a new $q^{k}$ by \eqref{eq:ewilq}, and then a new iterate of the field $u^{k+1}$ by \eqref{eq:ewifwd}. That is,
\begin{subequations}
  \begin{align}
    \lambda^{k} =& A^{-H}R^H(d-Ru^{k})\label{eq:recur1}\\
    q^k =& \frac{1}{\beta}(B^HB)^{-1}\lambda^k\\
     u^{k+1} =& A^{-1}q^k, \label{eq:recur2}
  \end{align}
\end{subequations}
which can be symbolically described in the following:
\begin{displaymath}
\cdots u^{k}  \overset{A^{-H}}{\rightarrow} \lambda^k \overset{(B^H B)^{-1}}{\rightarrow} q^k \overset{A^{-1}}{\rightarrow} u^{k+1} \overset{A^{-H}}{\rightarrow}\cdots
\end{displaymath}
By eliminating $q^k$, we have
\begin{equation}
  \begin{cases}
    \lambda^{k} = A^{-H}R^H(d-Ru^{k})\\
     u^{k+1} = \frac{1}{\beta}A^{-1}(B^HB)^{-1}\lambda^k
  \end{cases}
\end{equation}
or in matrix form
\begin{equation}
  \underbrace{\begin{bmatrix}
    I & 0\\
    \frac{1}{\beta}A^{-1}(B^HB)^{-1} & I
  \end{bmatrix}}_{M}\begin{bmatrix}
    \lambda^{k}\\
    u^{k+1}
  \end{bmatrix}=\underbrace{\begin{bmatrix}
  0 &-A^{-H} R^HR\\
  0 & 0
  \end{bmatrix}}_N\begin{bmatrix}
    \lambda^{k-1}\\
    u^k
  \end{bmatrix}+\begin{bmatrix}
  A^{-H} R^H d\\
  0
  \end{bmatrix}
\end{equation}
The diagonal and lower triangle of the matrix $M$, as well as the upper triangle of the matrix $N$ are none zero. The iterative scheme of this structure is exactly the well-known Gauss-Seidel iteration \citep[section 10.1.2]{Golub_1996_MATCOMP}. The Gauss-Seidel scheme may be costly but makes the use of both time and frequency domain solvers feasible for FWI using extended source.

In terms of \eqref{eq:recur1}, we plug in $\lambda^k = A^{-H} R^H(d-Ru^{k})$ into \eqref{eq:recur2}, yielding
\begin{equation}\label{eq:GS}
  Au^{k+1} = \frac{1}{\beta}(B^HB)^{-1}A^{-H}R^H(d-Ru^{k})
  = -\frac{1}{\beta}(B^HB)^{-1}A^{-H}R^H R u^{k} + \frac{1}{\beta}(B^HB)^{-1}A^{-H}R^H d
\end{equation}
Multiplying $A^{-1}$ on both side of \eqref{eq:GS}, the Gauss-Seidel recursion can be written as the fixed point iteration for single variable $u$
\begin{equation}
  u^{k+1}   = -\frac{1}{\beta}A^{-1}(B^HB)^{-1}A^{-H}R^H R u^{k} + \frac{1}{\beta}A^{-1}(B^HB)^{-1}A^{-H}R^H d. 
\end{equation}
Assume $u^*$ is the true solution after convergence, hence
\begin{equation}
  u^*   = -\frac{1}{\beta}A^{-1}(B^HB)^{-1}A^{-H}R^H R u^* + \frac{1}{\beta}A^{-1}(B^HB)^{-1}A^{-H}R^H d. 
\end{equation}
Define $e^k=u^k-u^*$ the error between current iterate with the true solution. It is then easy to establish the relation of the errors between two iterates
\begin{equation}
  \|e^{k+1}\| = \|u^{k+1}-u^*\|=\|\frac{1}{\beta}A^{-1}(B^HB)^{-1}A^{-H}R^H R(u^{k}-u^*)\|
  \leq \|\frac{1}{\beta}A^{-1}(B^HB)^{-1}A^{-H}R^H R\|\cdot\|e^k\|
\end{equation}
We conclude that regardless any initial value $u^0$, the Gauss-Seidel iterations are guaranteed to converge provided that the spectral radius satisfies
\begin{equation}\label{eq:spectral}
\|\frac{1}{\beta} A^{-1} (B^HB)^{-1} A^{-H} R^HR\| <1,
\end{equation}
This establishes a quantitative requirement to choose a proper penalty coefficient ($\beta>\|A^{-1} (B^HB)^{-1} A^{-H} R^HR\|$).

Assuming $A$ is invertible, for any positive $\beta$ the regularized Hessian $\beta (BA)^H (BA) + R^H R$ in \eqref{eq:normalewi} is always positive definite, therefore invertible for the solution of the linear system of $u$.  Unfortunately, the Gauss-Seidel scheme requires the penalty coefficient $\beta$ must be larger than certain values to converge.  It rules out the Gauss-Seidel iterations to seek for the solution of the misfit functional if $\beta$ is chosen very small, even though 
the solution of the resulting normal equation does exist.

\subsection{Recursion assisted by surrogate function}
 
To bypass the above limitation, we may introduce a surrogate function of the objective based on the wavefield at iterate $k$, inspired by the work of  \cite{aghamiry2020accurate}
\begin{equation}
  J'(u;u^k) = J(u) + \frac{1}{2}(u-u^k)^H (\gamma (BA)^H BA - R^H R) (u-u^k),
\end{equation}
where the parameter $\gamma$ has to be chosen such that $\gamma (BA)^H BA - R^H R$ is positive definite, thus $J'(u;u^k)\ge J(u)$ holds for any $u$.  Denote $u^{k+1}=\arg\min_u J'(u;u^k)$. It can be seen that
\begin{equation}
J(u^{k+1})\leq J'(u^{k+1}; u^k)\leq J'(u^k; u^k)=J(u^k).
\end{equation}
Consequently, the minimizers of the surrogate function at each iteration constructs a sequence of wavefield solution that approaches the minimizer of the original misfit function $J(u)$.

The surrogate function can be rewritten as
\begin{equation}
  \begin{split}
    J'(u;u^k)=& \frac{1}{2} \|d-Ru^k -R(u-u^k)\|^2 + \frac{\beta}{2}\|BA(u-u^k + u^k)\|^2 \\
    &+\frac{\gamma}{2}\|BA(u-u^k)\|^2 - \frac{1}{2} \|R (u-u^k)\|^2\\
  =&\frac{1}{2} \|d-Ru^k\|^2 -\langle d-Ru^k, R(u-u^k) \rangle + \frac{\beta+\gamma}{2}\|BA(u-u^k)\|^2 \\
  & +\frac{\beta}{2}\langle BA(u-u^k), BA u^k\rangle + \frac{\beta}{2}\|BAu^k\|^2
  \end{split}
\end{equation}
Minimization of the surrogate function $J'(u;u^k)$ requires $\partial_u J'=0$, yielding the normal equation
\begin{equation}
  \begin{split}
  -R^H(d-Ru^k) +(\beta+\gamma)A^H B^H BA(u-u^k) +\beta A^H B^H BA u^k =0\\
  \Leftrightarrow
  (\beta+\gamma) u=\gamma u^k +A^{-1} (B^H B)^{-1}A^{-H} R^H(d-Ru^k).
  \end{split}
\end{equation}

In terms of equations \eqref{eq:ewiadj}, \eqref{eq:ewilq} and \eqref{eq:newq2}, we can define the following intermediate fields to complete each step of the recursion
\begin{subequations}\label{eq:steps}
  \begin{align}
    \lambda^k :=& A^{-H} R^H(d-Ru^k)\\
    q^k := &\frac{1}{\beta}(B^H B)^{-1} \lambda^k\\
    \tilde{u}^k :=& A^{-1} q^k
    \end{align}
\end{subequations}
such that the minimizer of $J'(u;u^k)$ is given by
\begin{equation}\label{eq:ukp1}
u^{k+1} = \frac{1}{\beta+\gamma} (\beta\tilde{u}^k+\gamma u^k )
\end{equation}
The field $q$ can be stored using the same memory unit as $\lambda$ to avoid extra memory allocation. Each step of the above recursion requires at least two modellings (one adjoint modelling and one forward modelling). It is then symbolically described in the following:
\begin{displaymath}
  \cdots u^{k}  \overset{A^{-H}}{\rightarrow}
  \lambda^k  \overset{(B^HB)^{-1}}{\rightarrow}
  q^k \overset{A^{-1}}{\rightarrow}
  \tilde{u}^k\overset{\beta, \gamma, u^k}{\rightarrow} u^{k+1} \rightarrow\cdots
\end{displaymath}

The surrogate function is a quadratic form associated with operator $BA$ (its Hessian is $(BA)^H BA $). The convergence of the recursion will finally be dictated by the choice of $\gamma$ and the condition number
\begin{displaymath}
\text{cond}{(BA)} = \frac{\sigma_{\max}}{\sigma_{\min}},
\end{displaymath}
where $\sigma_{\max}$ and $\sigma_{\min}$ stand for maximum and minimum of the eigenvalues of $BA$. The choice of the annihilator $B$ will drastically affect the convergence rate. If $\gamma$ is too large, then $J'$ will significantly differ from $J$, the resulting convergence to the solution of $J$ would be arbitrarily slow. If $\gamma$ is too small, we have no guarantee that the surrogate function satisfies $J'(u;u^k)\geq J(u)$, hence the solution of surrogate function may not converge to the solution of $J(u)$.



\bibliographystyle{apalike}

\newcommand{\SortNoop}[1]{}

\end{document}